\input amstex
\magnification=1200
\documentstyle{amsppt}
\NoRunningHeads
\NoBlackBoxes
\topmatter
\title New mathematical methods for psychophysical filtering of experimental 
data\linebreak and their processing\endtitle
\author Denis V. Juriev\endauthor
\affil ul.Miklukho-Maklaya 20-180, Moscow 117437 Russia\linebreak
(e-mail: denis\@juriev.msk.ru)\endaffil
\date math.GM/0005275\enddate
\endtopmatter
\document This article, which is devoted to new mathematical methods for
psychophysical filtering of experimental data and their processing, may be
regarded as a development of ideas of Ref.1 on the experimental detection of
interactivity in complex systems. The main innovation is that a
psychophysical filtering is used consistently now for such detection.

\subhead 1. Basic procedure\endsubhead
Let us consider a natural, behavioral, social or economical system $\Cal S$.
It will be described by a set $\{\varphi\}$ of quntities, which characterize
it at any moment of time $t$ (so that $\varphi=\varphi_t$). One may suppose
that the evolution of the system is described by a differential equation 
$$\dot\varphi=\Phi(\varphi)$$
and look for the explicit form of the function $\Phi$ from the experimental
data on the system $\Cal S$. However, the function $\Phi$ may depend on time,
it means that there are some hidden parameters, which control the system
$\Cal S$ and its evolution is of the form
$$\dot\varphi=\Phi(\varphi,u),$$
where $u$ are such parameters of unknown nature. One may suspect that such 
parameters are chosen in a way to minimize some goal function $K$, which may 
be an integrodifferential functional of $\varphi_t$:
$$K=K(\left[\varphi_{\tau}\right]_{\tau\le t})$$
(such integrodifferential dependence will be briefly notated as 
$K=K([\varphi])$ below). More generally, the parameters $u$ may be divided
on parts $u=(u_1,\ldots,u_n)$ and each part $u_i$ has its own goal function
$K_i$. However, this hypothesis may be confirmed by the experiment very 
rarely. In the most cases the choice of parameters $u$ will seem accidental
or even random. Nevertheless, one may suspect that the controls $u_i$ are 
{\sl interactive}, it means that they are the couplings of the pure controls 
$u_i^\circ$ with the {\sl unknown or incompletely known\/} feedbacks:
$$u_i=u_i(u_i^\circ,[\varphi])$$
and each pure control has its own goal function $K_i$. Thus, it is
suspected that the system $\Cal S$ realizes an {\sl interactive game}.
There are several ways to define the pure controls $u_i^\circ$. One of them
was proposed in Ref.1. It is based on the integrodifferential filtration of
the controls $u_i$:
$$u^\circ_i=F_i([u_i],[\varphi]).$$
To verify the formulated hypothesis and to find the explicit form of the
convenient filtrations $F_i$ and goal functions $K_i$ one should use the
theory of interactive games, which supplies us by the predictions of the
game, and compare the predictions with the real history of the game for
any considered $F_i$ and $K_i$ and choose such filtrations and goal functions,
which describe the reality better. One may suspect that the dependence of
$u_i$ on $\varphi$ is purely differential for simplicity or to introduce the
so-called {\sl intention fields}, which allow to consider any interactive
game as differential. Moreover, one may suppose that
$$u_i=u_i(u_i^\circ,\varphi)$$
and apply the elaborated procedures of {\sl a posteriori\/} analysis and
predictions to the system.

In many cases this simple algorithm effectively unravels the hidden 
interactivity of a complex system, however, sometimes it does not work.
Therefore, more sophisticated procedures should be applied. One of them
will be described below.

Let us consider an interactive game with states $\psi$ and interactive 
controls $w=(w_1,\ldots, w_m)$ so that the evolution equation has the form
$$\dot\psi=\Psi(\psi,w).$$
The interactive controls are the couplings of pure controls $w^\circ=
(w^\circ_1,\ldots, w^\circ_m)$ with the unknown or incompletely known 
feedbacks:
$$w=w(w^\circ,\psi).$$
Often such interactive game is a game of a real human person and the 
interactivity is caused by the coupling of conscious and subconscious
behavioral reactions.

Let us consider the integrodifferential nonlinear operator, which will
expresses the pure controls $w^\circ$ via $u$ and $\varphi$:
$$w^\circ=W([u],[\varphi]).$$
In the simplest case one may put
$$w^\circ=W(u,\varphi).$$
The pure controls $w^\circ$ realize a scenario for the interactive game
and it becomes a performance.

The integrodifferential filtration should be applied to $w$ instead of
$u$ now:
$$u^\circ_i=F_i([w],[\varphi],[\psi]).$$
We suspect that $u^\circ=(u^\circ_1,\ldots,u^\circ_n)$ are pure controls
for $u=(u_1,\ldots,u_n)$. It means that 
$$u=u(u^\circ,\varphi)$$
and the explicit dependence is unknown or incompletely known.
Because we use a human person, the whole procedure realizes 
{\it a psychophysical filtering\/} of the experimental data. 

\subhead 2. Psychophysical engineering\endsubhead

The basic procedure, which is described above, admits further improvements,
which are based on the ideas of automata theory and form the {\it 
psychophysical engineering}.

First, let us mark that an interactive game may be regarded as an infinite
automata with continuous time. Therefore, one is able to consider a
{\sl composition\/} of interactive games. It means that the interactive
controls of one game are used as pure controls for other one. Certainly,
some integrodifferential transforms of controls (analogous to one decribed
above) may be used between output of one game and the input of another to 
glue the games. Such transforms are defined by the {\sl known\/} 
integrodifferential nonlinear operators. Of course, such operators are 
reduced to functions in the simplest cases. 

Such glueings of interactive games allows to construct very sophisticated
objects started from only few basic interactive games. Combined with the
basic procedure above such psychophysical engineering becomes an effective
tool for the unraveling of interactivity in complex systems. Concrete parts
of the whole construction may be realized either by one human person or
by their group. Of course, external interactive phenomena of non-human nature
may be used also. 

The simplest interactive games, which may be used for the psychophysical
engineering, are the visual (or general perceptive) games. However, it is 
difficult to parallelize such games so it is reasonable to {\sl interiorize\/}
the behavioral reactions, which are generated in such games. It means that the
functional dependences between pure and interactive controls are remebered and 
reproduced by the person without any external or even internal visual image.

Of course, the process of psychophysical engineering may be regarded as
a {\sl tactical procedure\/} [2] and may contain the intuitive subconscious 
steps. Therefore, a psychophysiological training may sufficiently enlarge the
capabilities of a person to perform the psychophysical engineering. 
The group trainings are effective when groups are used. Personal abilities 
to adapt the external interactivity of non-human nature are also very 
essential.

\Refs
\roster
\item"[1]" Juriev D., Experimental detection of interactive phenomena and
their analysis. E-print: math.GM/0003001.
\item"[2]" Juriev D., Tactical games \& behavioral self-organization.
E-print: math.GM/9911147.
\endroster
\endRefs
\enddocument